# Non-Degenerate Conditionings of the Exit Measures of Super Brownian Motion

By
Thomas S. Salisbury and John Verzani
*York University and the Fields Institute*
*CUNY – College of Staten Island*

*Summary:*
We introduce several martingale changes of measure of the law of the exit measure of super Brownian motion. These changes of measure include and generalize one arising by conditioning the exit measures to charge a point on the boundary of a 2-dimensional domain. In the case we discuss this is a non-degenerate conditioning. We give characterizations of the new processes in terms of "immortal particle" branching processes with immigration of mass, and give applications to the study of solutions to $Lu = cu^2$ in $D$. The representations are related to those in an earlier paper, which treated the case of degenerate conditionings.

## 1. INTRODUCTION

We investigate conditionings of the exit measures of super Brownian motionin $\mathbb{R}^d$. We can think of super Brownian motion as the limit of a particle system, which can heuristically be described as follows. It consists of a cloud of particles, each diffusing as a Brownian motion and undergoing critical branching. A measure valued process is formed by assigning a small point mass to each particle's position at a given time. The exit measure $X^D$ from a domain $D$ is then formed by freezing this mass at the point the particle first exits from $D$. For a sequence of subdomains, these measures can be defined on the same probability space, giving rise to a process indexed by the subdomains. In dimension 2, with positive probability, points on the boundary of a smooth enough domain will be hit by the support of the exit measure. In this paper, we study conditionings of the sequence of exit measures, analogous to the conditioning by this event. Unlike the case $d = 2$, in higher dimensions the corresponding event has probability 0, and the analogous conditioning is a degenerate one. Such degenerate conditionings were treated in the paper [12] (SV1).

To be more specific, let $D$ be a bounded domain in dimension $d = 2$, and let $D_k$ be an increasing sequence of subdomains. The domains $D_k$ give rise to a process of exit measures $X^k$, each defined on the boundary of $D_k$. We work under $\mathbb{N}_x$, the excursion measure under which Le Gall's Brownian snake evolves. Let $\hat{\mathbb{M}}_x$ be the law of super Brownian motion, conditioned on the exit measure hitting a fixed point $z$ on $\partial D$ (that is, conditioned on it charging all balls containing $z$). Let $\mathcal{F}_k$ be the $\sigma$-field generated by the particles before they exit $D_k$ and denote integration by $\langle,\rangle$. Our first result is an explicit description of $\hat{\mathbb{M}}_x$ on $\mathcal{F}_k$. Its densities with respect to $\mathbb{N}_x$ form a martingale (in $k$) which can be explicitly written in terms of the $X^k$.

More generally, the differential equation $Lu = 4u^2$ plays an important role in our discussion, and for the exit measures in general. In Lemma 3.1 it is shown that if $g \geq u \geq 0$ are both solutions in $D$ to $Lu = 4u^2$ then $\hat{M}_k = \exp-\langle X^k, u\rangle - \exp-\langle X^k, g\rangle$ is an $\mathcal{F}_k$ martingale. Letting $v = g - u$, we can define a general change of measure, using this martingale, to give a measure $\hat{\mathbb{M}}_x$ satisfying

$$\frac{d\hat{\mathbb{M}}_x}{d\mathbb{N}_x}\bigg|_{\mathcal{F}_k} = \frac{1}{v(x)}\hat{M}_k$$

for each $k$. (In the above example, $u = 0$ and $g = g_z = \mathbb{N}_x(\mathcal{R}^D \cap \{z\} \in \partial D)$, where $\mathcal{R}^D$ denotes the range of the super Brownian "particles" before exiting $D$.)

Our second result is that the measures $\hat{\mathbb{M}}_x$ on $\mathcal{F}_k$ can be represented in terms of a branching process of "immortal particles" together with immigration of mass. Two such equivalent representations are given, in Theorems 3.2 and 3.5. The first involves a conditioned diffusion in which particles may die, but when this occurs two independent particles are born as replacements. The other uses a conservative conditional diffusion undergoing binary branching. The branching mechanism in both representations is homogeneous, unlike the representations of SV1. (We use the terminology "immortal particle" to refer to a particle that is conditioned to exit $D$ through $\partial D$. The language comes from previous conditionings of the superprocess. See SV1 for references to earlier work.)

By using both descriptions, we can investigate the solutions to the equation $Lu = 4u^2$. We see in two examples that the solutions given by

$$g_z(x) = \mathbb{N}_x(\mathcal{R}^D \cap \{z\})$$

for $z \in \partial D$ and by

$$g_f(x) = \mathbb{N}_x(1 - \exp(-\langle X^d, f\rangle))$$

lead to quite different immortal particles pictures: the former having infinitely many branches and the latter just finitely many (when $f$ is bounded). We conjecture that the class of moderate functions (as studied in [10] and [6]), that is those solutions bounded in $D$ by a harmonic function, are precisely those for which the immortal particle picture has finitely many branches.

Finally, we draw an analogy between these conditionings and those treated in SV1. (See remark 3.12 for a description of the results in SV1.) In that paper we investigated transforms based on a different type of martingale than used here. That family of martingales generalized the ones arising from conditioning the exit measure to hit a given finite number of points on the boundary of $D$, in the case that this conditioning was degenerate (that is, that the event conditioned on had probability 0). Because of this degeneracy, the results there had an asymptotic character, and required analytic estimates of small solutions to certain nonlinear PDEs. Those conditionings also had immortal particle representations, though the particles in their backbones evolved in an inhomogeneous manner. In section 4 of the current paper, we present a martingale change of measure combining features of the conditionings $\hat{\mathbb{M}}_x$ described above, and those of SV1. In Theorem 4.4 we derive an immortal particle representation for this general class of transforms.

2. PRELIMINARIES

This paper is a sequel to SV1 [12], but for the convenience of the reader we restate some of the lemmas used therein. Any proofs appear in SV1.

2.1. **Notation.** For a set $A$, let $|A|$ denote its cardinality, and let $\mathcal{P}(A)$ denote the collection of partitions of $A$. Choose some arbitrary linear order $\prec$ on the set of finite subsets of the integers. For $A$ such a finite subset, and $\sigma \in \mathcal{P}(A)$, let $\sigma(j)$ be the $j$th element of $\sigma$ in this order. Thus for example,

$$\prod_{C \in \sigma} \langle X^k, v^C \rangle = \prod_{j=1}^{|\sigma|} \langle X^k, v^{\sigma(j)} \rangle.$$

We will switch between these notations according to which seems clearer.

2.2. **Set facts.** We make use of the following two simple lemmas. We will use the convention that a sum over an empty set is 0.

**Lemma 2.1** (Lemma 2.1 of SV1). *Let $A \subseteq B \subseteq C$ be subsets of $\{1, 2, \ldots n\}$. Then*

$$\sum_{A \subseteq B \subseteq C} (-1)^{|B|} = (-1)^{|C|} \mathbb{1}_{A=C}.$$

**Lemma 2.2** (Lemma 2.2 of SV1). *Let $A$ be finite, and let $w_i \in \mathbb{R}$ for $i \in A$. Then*

$$\prod_{i \in A}(1 - w_i) = 1 + \sum_{\substack{C \subseteq A \\ \emptyset \neq C}} (-1)^{|C|} \Big( \prod_{i \in C} w_i \Big).$$

In this paper we use the letter $K$ to denote a generic non-trivial constant whose particular value may vary from line to line. If it is important, explicit dependencies on other values will be specified.

2.3. **Facts about conditioned diffusions.** First we recall some formulae for conditioned Brownian motion.

Let $B$ be $d$-dimensional Brownian motion started from $x$, under a probability measure $P_x$. Write $\tau_D$ for the first exit time of $B$ from $D$.

Let $g : D \to [0, \infty)$ be bounded on compact subsets of $D$, and set

$$L_g = \frac{1}{2}\Delta - g.$$

Let $\xi_t$ be a process which, under a probability law $P_x^g$, has the law of a diffusion with generator $L_g$ started at $x$ and killed upon leaving $D$. In other words, $\xi$ is a Brownian motion on $D$, killed at rate $g$. Write $\zeta$ for the lifetime of $\xi$. Then

$$E_x^g(\xi_t \in A, \zeta > t) = E_x(\exp - \int_0^t ds\, g(B_s), B_t \in A, \tau_D > t). \qquad (2.1)$$

Let $U^g f(x) = \int_0^\infty P_x^g(f(\xi_t) \mathbf{1}_{\{\zeta > t\}}) dt$ be the potential operator for $L_g$. If $g = 0$ we write $U$ for $U^g$. If $0 \leq u$ is $L_g$-superharmonic, then the law of the $u$-transform of $\xi$ is determined by the formula

$$P_x^{g,u}(\Phi(\xi)\mathbf{1}_{\{\zeta > t\}}) = \frac{1}{u(x)} P_x^g(\Phi(\xi) u(\xi_t) \mathbf{1}_{\{\zeta > t\}})$$

for $\Phi(\xi) \in \sigma\{\xi_s; s \leq t\}$. Assuming that $0 < u < \infty$ on $D$, this defines a diffusion on $D$. If $u$ is $L_g$-harmonic, then it dies only upon reaching $\partial D$. If for $f \geq 0$, $u = U^g f$ (that is, $u$ is a potential) then it dies in the interior of $D$, and in fact $P_x^{g,u}$ satisfies

$$P_x^{g,u}(\Phi(\xi)) = \frac{1}{u(x)} \int_0^\infty P_x^g(\Phi(\xi_{\leq t}) f(\xi_t) \mathbf{1}_{\{\zeta > t\}}) dt, \tag{2.2}$$

where $\xi_{\leq t}$ is the process $\xi$ killed at time $t$.

In addition, if $u = h + v$, where $h \geq 0$ is $L_g$-harmonic, and $v = U^g f$ with $f \geq 0$, then

$$P_x^{g,u} = \frac{1}{u(x)} \left( h(x) P_x^{g,h} + v(x) P_x^{g,v} \right). \tag{2.3}$$

Suppose that $\Phi$ vanishes on paths $\xi$ that reach $\partial D$. Applying (2.2) to (2.3), we see that

$$P_x^{g,u}(\Phi(\xi)) = \frac{v(x)}{u(x)} P_x^{g,v}(\Phi(\xi)) = \frac{1}{u(x)} \int_0^\infty P_x^g(\Phi(\xi_{\leq t}) f(\xi_t) \mathbf{1}_{\{\zeta > t\}}) dt, \tag{2.4}$$

in this case as well.

**2.4. Facts about the Brownian snake.** Next we recall some useful facts about the Brownian snake. Refer to [1] or [3] for a general introduction to superprocesses.

The Brownian snake is a path-valued process, devised by Le Gall as a means to construct super Brownian motion without limiting procedures. Refer to [7] or [9] for the construction.

We use the standard notation $(W_s, \zeta_s)$ for the Brownian snake, and $\mathbb{N}_x$ for the excursion measure of the Brownian snake starting from the trivial path $(w, \zeta), \zeta = 0, w(0) = x$. Note that $W_s(\cdot)$ is constant on $[\zeta_s, \infty)$, and $\zeta$ has the distribution of a Brownian excursion under $\mathbb{N}_x$.

Super Brownian motion $X_t$ is defined as

$$\langle X_t, \phi \rangle = \int \phi(W_s(t)) \, dL_t(s),$$

where $L_t$ is the local time of $\zeta$ at level $t$. Dynkin [4] introduced the exit measure $X^D$ associated with $X_t$. We follow Le Gall's snake-based definition of $X^D$ (see [9]) as

$$\langle X^D, \phi \rangle = \int \phi(W_s(\zeta_s)) \, dL^D(s),$$

where $L^D(\cdot)$ is an appropriate local time for $W_s(\zeta_s)$ on $\partial D$.

We denote the range of the Brownian snake by $\mathcal{R}(W) = \{W_s(t) : 0 \leq s \leq \sigma, 0 \leq t \leq \zeta_s\}$ and the range inside $D$ by $\mathcal{R}^D(W) = \{W_s(t) : 0 \leq s \leq \sigma, 0 \leq t \leq \tau_D(W_s) \wedge \zeta_s\}$. There is an obvious inclusion between the range inside $D$ and the exit measures, given by

$$\{\langle X^D, \mathbf{1}_A \rangle > 0\} \subseteq \{\mathcal{R}^D(W) \cap A \neq \emptyset\}.$$

We refer the reader to Le Gall [9] for other facts about the Brownian snake, including the following result (an immediate consequence of Theorem 4.2 and its corollary in [9]).

**Lemma 2.3.** *Let $g$ be a solution to $\Delta g = 4g^2$ in $D$, and let $\{D_k\}$ be an increasing sequence of smooth subdomains of $D$. Then for each $k$,*

$$\mathbb{N}_x(1 - \exp-\langle X^{D_k}, g\rangle) = g(x).$$

Let $\mathcal{F}_k = \mathcal{F}_{D_k}$ be the $\sigma$-field of events determined by the superprocess killed upon exiting $D_k$. See [4] for a formal definition. Or refer to the final section of SV1, which gives a definition in terms of the historical superprocess.

Dynkin introduced a Markov property for the exit measures in [4]. In our context, the Markov property is established in [10]. The next result gives it in the form we will use it:

**Lemma 2.4.**

$$\mathbb{N}_x(\exp-\langle X^D, \phi\rangle \mid \mathcal{F}_k) = \exp-\langle X^{D_k}, \mathbb{N}.(1 - \exp-\langle X^D, \phi\rangle)\rangle.$$

We use the following notation, where $B_s$ denotes a path in $D$ whose definition will be clear from the context:

$$e^D_\phi = e^D(\phi) = \exp-\langle X^D, \phi\rangle,$$
$$\mathcal{N}_t(e^D_\phi) = \mathcal{N}_t(e^D_\phi, B) = \exp-\int_0^t ds\, 4\mathbb{N}_{B_s}(1 - e^D_\phi).$$

The Palm formula for the Brownian snake takes the form: (cf. [9], Proposition 4.1)

$$\mathbb{N}_x(\langle X^D, \phi\rangle e^D_\psi) = E_x(\phi(B_{\tau_D})\mathcal{N}_{\tau_D}(e^D_\psi)). \tag{2.5}$$

We will make use of the following extension to the basic Palm formula. See [2] for a general discussion of this type of Palm formula.

**Lemma 2.5** (Lemma 2.5 of SV1). *Let $N = \{1, 2, \ldots n\}$, $n \geq 2$. Let $D$ be a domain, and let $B$ be a Brownian motion in $D$ with exit time $\tau$. Let $\{\psi_i\}$ be a family of measurable functions. Then*

$$\mathbb{N}_x\Big(e_\phi \prod_{i\in N}\langle X^D, \psi_i\rangle\Big)$$
$$= \frac{1}{2}\sum_{\substack{M\subseteq N \\ \emptyset, N\neq M}} E_x\Big(4\int_0^\tau dt\, \mathcal{N}_t(e_\phi)\mathbb{N}_{B_t}\Big(e_\phi\prod_{i\in M}\langle X^D, \psi_i\rangle\Big)\mathbb{N}_{B_t}\Big(e_\phi\prod_{i\in N\setminus M}\langle X^D, \psi_i\rangle\Big)\Big).$$

Using the extended Palm formula one may show an exponential bound on the moments of the exit measure.

**Lemma 2.6** (Lemma 2.7 of SV1). *Let $D$ be a domain in $\mathbb{R}^d$ satisfying $\sup_{x\in D} E_x(\tau_D) < \infty$, where $\tau_D$ is the exit time from $D$ for Brownian motion. Then there exists $\lambda > 0$ such that*

$$\sup_{x\in D} \mathbb{N}_x(\exp\lambda\langle X^D, 1\rangle - 1) < \infty.$$

*Remark* 2.7. A bounded domain $D$ in $\mathbb{R}^d$ will satisfy $\sup_D E_x(\tau_D) < \infty$.

## 3. The martingale $\hat{M}_k$.

In this section we investigate a type of $h$-transform of the exit measures, given by a martingale change of measure. This transform is then interpreted in terms of a branching system of particles, as in SV1, but unlike the situation there, the branching system is now a homogeneous one.

Suppose $D$ is a bounded domain in $\mathbb{R}^d$, $D_k$ are smooth domains satisfying $D_k \Uparrow D$, and $g \geq 0$ satisfies $\frac{1}{2}\Delta g = 2g^2$ in $D$. Let $u \geq 0$ be a second solution to this equation, with $u \leq g$. Set $v = g - u$, and let $\hat{M}_k = \exp -\langle X^k, u \rangle - \exp -\langle X^k, g \rangle$.

**Lemma 3.1.** $\hat{M}_k$ is a $\mathcal{F}_k$ martingale.

*Proof.* Let $j < k$. Then by Lemmas 2.3 and 2.4,

$$\begin{aligned}
\mathbb{N}_x(\hat{M}_k \mid \mathcal{F}_j) &= \mathbb{N}_x(\exp -\langle X^k, u \rangle - \exp -\langle X^k, g \rangle \mid \mathcal{F}_j) \\
&= \exp -\langle X^j, \mathbb{N}.(1 - \exp -\langle X^k, u \rangle)\rangle - \exp -\langle X^j, \mathbb{N}.(1 - \exp -\langle X^k, g \rangle)\rangle \\
&= \exp -\langle X^j, u \rangle - \exp -\langle X^j, g \rangle = \hat{M}_j.
\end{aligned}$$

□

As a consequence, we can define a transformed process via a martingale change of measure. If $\Phi_k$ is a $\mathcal{F}_k$-measurable function, set

$$\hat{\mathbb{M}}_x(\Phi_k) = \frac{1}{v(x)} \mathbb{N}_x(\Phi_k \hat{M}_k) = \frac{1}{v(x)} \mathbb{N}_x(\Phi_k e_g^k (e^{\langle X^k, v \rangle} - 1))$$

$$= \frac{1}{v(x)} \mathbb{N}_x(\Phi_k e_g^k \sum_{n=1}^\infty \frac{1}{n!} \langle X^k, v \rangle^n). \tag{3.1}$$

As in SV1 we can represent these conditioned exit measures in terms of a branching process with immigration.

First we construct a homogeneous branching process. Our underlying process will be Brownian motion killed at rate $4g$ (that is, with generator $L_{4g}$). Recall that $E^{4g}$ denotes its law. Then $v$ satisfies

$$L_{4g} v = \frac{1}{2}(\Delta g - \Delta u) - 4gv = 2v(g + u - 2g) = -2v^2.$$

In other words, $v$ is $L_{4g}$-superharmonic, and we can consider the $v$-transform of the $L_{4g}$-process. Recall that its law is denoted by $E^{4g,v}$.

To describe the branching $g$-process, we start with a single particle, with law $E_x^{4g,v}$. When it dies, say at $y$, it is replaced by two independent offspring with laws $E_y^{4g,v}$. In other words, those offspring will now evolve with the same transition function as their parent, as will all their descendants. Again denote the branching process viewed as a measure on $D$ by $\Upsilon$, and let $\Upsilon^k$ be the measure generated by those particle which haven't left $D_k$. Let $\mathbb{Q}_x$ denote the law of $\Upsilon$.

It is worth pointing out that if $v$ is an $L_{4g}$-potential, then the above branching process will have infinitely many branches, as no particle ever makes it out to $\partial D$. If $v$ has an $L_{4g}$-harmonic component, then some particles may exit $D$, and the total number of branches will be finite with positive probability.

We immigrate mass along the branches of the backbone traced out by the branching particle system to construct exit measures $Y^k$ on $D_k$ for each $k$. This yields a measure $\hat{\mathbb{N}}_x$, under which $Y^k$ has law defined by

$$\hat{\mathbb{N}}_x(\exp-\langle Y^k, \phi\rangle) = \mathbb{Q}_x(\exp - \int_0^\infty dt\, 4\langle \Upsilon_t^k, \tilde{\mathbb{N}}.(1-e_\phi^k)\rangle),$$

where $\tilde{\mathbb{N}}_y(\Phi_k) = \mathbb{N}_y(\Phi_k \exp-\langle X^k, g\rangle)$. In other words, given the backbone $\Upsilon^k$, we form a Poisson random measure $N^k(d\mu)$ with intensity $\int_0^\infty dt \int 4\Upsilon_t^k(dy)\tilde{\mathbb{N}}_y(X^k \in d\mu)$. We then realize the exit measure under $\hat{\mathbb{N}}_x$ as $Y^k = \int \mu N^k(d\mu)$. As in SV1, $\tilde{\mathbb{N}}_y$ is also the excursion measure for the snake based on Brownian motion killed at rate $g$.

Since the branching process is homogeneous, we can partition the particles into classes determined by their having a common ancestor prior to exiting $D_k$. Let $\Upsilon^k \sim n$ denote the event that there are $n$ distinct ancestors before exiting $D_k$.

**Theorem 3.2.** *We have, in the notation of this section,*

$$\hat{\mathbb{M}}_x(\exp-\langle X^k, \phi\rangle) = \hat{\mathbb{N}}_x(\exp-\langle Y^k, \phi\rangle).$$

*Remark* 3.3. Using historical processes, as in the last section of SV1 one can show that $\hat{\mathbb{M}}_x = \hat{\mathbb{N}}_x$ on $\mathcal{F}_k$.

*Proof.* We show by induction that

$$v(x)\mathbb{Q}_x\left(\exp - \int_0^\infty dt\, 4\langle \Upsilon_t^k, \tilde{\mathbb{N}}.(1-e_\phi^k)\rangle, \Upsilon^k \sim n\right) = \frac{1}{n!}\mathbb{N}_x(e_{\phi+g}^k\langle X^k, v\rangle^n). \qquad (3.2)$$

From this it follows, by summing on $n$, that

$$\hat{\mathbb{N}}_x(\exp-\langle Y^k, \phi\rangle) = \frac{1}{v(x)}\sum_{n=1}^\infty \frac{1}{n!}\mathbb{N}_x(e_{\phi+g}^k\langle X^k, v\rangle^n)$$

$$= \mathbb{N}_x(e_\phi^k \exp-\langle X^k, g\rangle\left(\exp(\langle X^k, v\rangle) - 1\right)) = \hat{\mathbb{M}}_x(e_\phi^k).$$

This will prove the theorem.

We note that by Lemma 2.5, each term above satisfies

$$\frac{1}{n!}\mathbb{N}_x(e_{\phi+g}^k\langle X^k, v\rangle^n)$$

$$= \frac{1}{n!}\sum_{j=1}^{n-1}\binom{n}{j}E_x\left(2\int_0^{\tau_k}\mathcal{N}_t(e_{\phi+g}^k)\mathbb{N}_{B_t}(e_{\phi+g}^k\langle X^k, v\rangle^j)\mathbb{N}_{B_t}(e_{\phi+g}^k\langle X^k, v\rangle^{n-j})dt\right)$$

$$= \sum_{j=1}^{n-1}E_x(2\int_0^{\tau_k}\mathcal{N}_t(e_{\phi+g}^k)\left(\frac{1}{j!}\mathbb{N}_{B_t}(e_{\phi+g}^k\langle X^k, v\rangle^j)\right)\left(\frac{1}{(n-j)!}\mathbb{N}_{B_t}(e_{\phi+g}^k\langle X^k, v\rangle^{n-j})\right)dt). \qquad (3.3)$$

Now to establish (3.2). First, in the case when $n=1$ we have that $\Upsilon^k$ is given by a single $v$-process which has lifetime greater than $\tau_k$. Hence,

$$v(x)\mathbb{Q}_x(\exp - \int_0^\infty dt\, 4\langle \Upsilon_t^k, \tilde{\mathbb{N}}.(1-e_\phi^k)\rangle, \Upsilon^k \sim 1)$$

$$= v(x)E_x^{4g,v}(\exp - \int_0^{\tau_k} dt\, 4\tilde{\mathbb{N}}_{\xi(t)}(1-e_\phi^k), \zeta > \tau_k)$$

$$
\begin{aligned}
&= E_x^{4g}\left(v(\xi(\tau_k))\exp-\int_0^{\tau_k} dt\, 4\tilde{\mathbb{N}}_{\xi(t)}(1-e_\phi^k), \zeta>\tau_k\right)\\
&= E_x(v(\xi(\tau_k))\exp\left(-\int_0^{\tau_k} dt\, 4g(\xi(t))\right)\exp\left(-\int_0^{\tau_k} dt\, 4\mathbb{N}_{\xi(t)}(e_g^k(1-e_\phi^k))\right))\\
&= E_x(v(\xi(\tau_k))\exp-\int_0^{\tau_k} dt\, 4\Big(g(\xi(t))-\mathbb{N}_{\xi(t)}(1-e_g^k)+\mathbb{N}_{\xi(t)}(1-e_{\phi+g}^k)\Big))\\
&= E_x(v(\xi(\tau_k))\mathcal{N}_{\tau_k}(e_{\phi+g}^k)) \tag{3.4}\\
&= \mathbb{N}_x(e_{\phi+g}^k\langle X^k, v\rangle), \tag{3.5}
\end{aligned}
$$

where (3.4) follows from (2.3) and (3.5) from (2.5).

When $n>1$ the first particle splits at its lifetime $\zeta<\tau_k$. By the Markov property for $\Upsilon$ and the conditional independence of the offspring,

$$
\begin{aligned}
&v(x)\mathbb{Q}_x(\exp-\int_0^\infty dt\, 4\langle \Upsilon_t^k, \tilde{\mathbb{N}}.(1-e_\phi^k)\rangle, \Upsilon^k\sim n)\\
&=\sum_{j=1}^{n-1} v(x) E_x^{4g,v}(\exp\left(-\int_0^\zeta dt\, 4\tilde{\mathbb{N}}_{\xi(\zeta)}(1-e_\phi^k)\right)\\
&\qquad\times \mathbb{Q}_{\xi(\zeta)}\left(\exp-\int_0^\infty dt\, 4\langle \Upsilon_t^k, \tilde{\mathbb{N}}.(1-e_\phi^k)\rangle, \Upsilon^k\sim n-j\right)\\
&\qquad\times \mathbb{Q}_{\xi(\zeta)}\left(\exp-\int_0^\infty dt\, 4\langle \Upsilon_t^k, \tilde{\mathbb{N}}.(1-e_\phi^k)\rangle, \Upsilon^k\sim j\right), \zeta<\tau_k)\\
&=\sum_{j=1}^{n-1} E_x^{4g}(\int_0^{\tau_k} ds\, 2v^2(\xi(s))\mathbb{1}_{\zeta>s}\exp\left(-\int_0^s dt\, 4(\tilde{\mathbb{N}}_{\xi(\zeta)}(1-e_{\phi+g}^k)-(\tilde{\mathbb{N}}_{\xi(\zeta)}(1-e_g^k))\right)\\
&\qquad\times\mathbb{Q}_{\xi(s)}\left(\exp-\int_0^\infty dt\, 4\langle\Upsilon_t^k,\tilde{\mathbb{N}}.(1-e_\phi^k)\rangle, \Upsilon^k\sim n-j\right)\\
&\qquad\times\mathbb{Q}_{\xi(s)}\left(\exp-\int_0^\infty dt\, 4\langle\Upsilon_t^k,\tilde{\mathbb{N}}.(1-e_\phi^k)\rangle, \Upsilon^k\sim j\right))
\end{aligned}
\tag{3.6}
$$

$$
\begin{aligned}
&=\sum_{j=1}^{n-1} E_x(2\int_0^{\tau_k} ds\,\exp\left(-\int_0^s dt\, 4g(\xi(t))\right)\exp\left(\int_0^s dt\, 4g(\xi(t))\right)\mathcal{N}_s(e_{\phi+g}^k)\\
&\qquad\times v(\xi(s))\mathbb{Q}_{\xi(s)}\left(\exp-\int_0^\infty dt\, 4\langle\Upsilon_t^k,\tilde{\mathbb{N}}.(1-e_\phi^k)\rangle, \Upsilon^k\sim n-j\right)\\
&\qquad\times v(\xi(s))\mathbb{Q}_{\xi(s)}\left(\exp-\int_0^\infty dt\, 4\langle\Upsilon_t^k,\tilde{\mathbb{N}}.(1-e_\phi^k)\rangle, \Upsilon^k\sim j\right))
\end{aligned}
\tag{3.7}
$$

$$
\begin{aligned}
&=\sum_{j=1}^{n-1} E_x(2\int_0^{\tau_k} ds\,\mathcal{N}_s(e_{\phi+g}^k)\left(\frac{1}{(n-j)!}\mathbb{N}_{\xi(s)}(e_{\phi+g}^k\langle X^k,v\rangle^{n-j})\right)\\
&\qquad\times\left(\frac{1}{j!}\mathbb{N}_{\xi(s)}(e_{\phi+g}^k\langle X^k,v\rangle^j)\right))
\end{aligned}
\tag{3.8}
$$

$$
=\frac{1}{n!}\mathbb{N}_x(e_{\phi+g}^k\langle X^k,v\rangle^n) \tag{3.9}
$$

Line (3.6) follows from (2.4), (3.7) follows from (2.1) and Lemma 2.3, (3.8) by the inductive hypotheses, and (3.9) by (3.3). □

*Remark* 3.4. The relation $L_{4g}v = -2v^2$ is closely related to the following recurrence, derived from equation (3.2) of SV1.

$$L_{4g}v^A = -2 \sum_{\substack{B \subseteq A \\ \emptyset, A \neq B}} v^B v^{A \setminus B}. \tag{3.10}$$

To see this, fix $n$ and let $v^A = c_{|A|}v$. Then (3.10) holds for $|A| > 1$ (though not for $|A| = 1$), provided

$$c_k = \sum_{j=1}^{k-1} \binom{k}{j} c_j c_{k-j} \tag{3.11}$$

for $1 < k \leq n$. As in the proof of Lemma 3.1 of Serlet [13] (or of Lemma 2.7 of SV1), the latter has solution

$$c_{k+1} = a^{k+1} \frac{2^k (2k)!}{k!},$$

where $a$ is arbitrary. In this case, the $M^k$ of (3.21) becomes

$$e_g^k \sum_{\sigma \in \mathcal{P}(N)} \prod_{A \in \sigma} \langle X^k, v^A \rangle = e_g^k \sum_{m=1}^{n} \langle X^k, v \rangle^m \sum_{\substack{\sigma \in \mathcal{P}(N) \\ |\sigma| = m}} \prod_{A \in \sigma} c_{|A|}.$$

Using (3.11) repeatedly, this is easily seen to equal

$$e_g^k \sum_{m=1}^{n} \frac{c_N}{m!} \langle X^k, v \rangle^m.$$

Choosing $a$ to make $c_N = 1$ produces a truncated version of the martingale $\hat{M}^k$.

There is an alternative description of the above backbone, which is in some ways more natural, though it is less closely tied to the approach of SV1. In this version, the backbone is a branching diffusion. The diffusion is again a $v$-transform, but this time of the process with generator $L_{2(u+g)}$. Note that now $L_{2(u+g)}v = 0$. We denote this process by $\lambda_t$. We then let $\lambda$ branch at rate $2v$, to produce a tree $\Upsilon_t$. We write $\check{Q}_x$ for its law. The branching may be thought of as killing at rate $2v$, except that at the death time $\zeta$, instead of dying we branch in two. Branching is of course done independently. On top of this branching process, we immigrate mass exactly as before, to produce a measure $\check{\mathbb{N}}_x$. It turns out to be the same as the measure $\hat{\mathbb{N}}_x$ given above.

**Theorem 3.5.** *For the measures described above one has*

$$\check{\mathbb{N}}_x(e_\phi^k) = \hat{\mathbb{M}}_x(e_\phi^k).$$

*Remark* 3.6. Actually we will show that $\check{\mathbb{N}}_x = \hat{\mathbb{N}}_x$, so that it follows from Remark 3.3 that $\check{\mathbb{N}}_x = \hat{\mathbb{M}}_x$ on $\mathcal{F}_k$.

*Proof.* The above statement could be proved directly, in a manner similar to the proof of Theorem 3.2. Instead, we will prove a slightly stronger statement, namely that $\Upsilon^k$ has the same law under $\hat{Q}_x$ as under $\mathbb{Q}_x$, for every $k$. The conclusion of the Theorem will therefore follow from Theorem 3.2.

To see this, it will suffice to show that, for $\phi_t(x) \geq 0$ and measurable in $(t, x)$, we have

$$\hat{Q}_x(\exp - \int_0^\infty dt \, \langle \Upsilon_t^k, \phi_t \rangle) = \mathbb{Q}_x(\exp - \int_0^\infty dt \, \langle \Upsilon_t^k, \phi_t \rangle). \tag{3.12}$$

For simplicity, write $A_x(n) = \hat{Q}_x(\exp - \int_0^\infty dt \, \langle \Upsilon_t^k, \phi_t \rangle; \Upsilon^k \sim n)$. We will show, by induction, that

$$A_x(n) = \mathbb{Q}_x(\exp - \int_0^\infty dt \, \langle \Upsilon_t^k, \phi_t \rangle; \Upsilon^k \sim n).$$

The result then follows by summing on $n$.

First the case $n = 1$:

$$A_x(1) = \hat{Q}_x \left( \exp - \int_0^\infty dt \, \langle \Upsilon_t^k, \phi_t \rangle; \Upsilon^k \sim 1 \right)$$

$$= E_x^{2(u+g),v} \left( \exp(- \int_0^{\tau_k} ds \, \phi_s(\lambda_s)); \zeta > \tau_k \right)$$

$$= E_x^{2(u+g),v} \left( \exp(- \int_0^{\tau_k} ds \, 2v(\lambda_s)) \exp(- \int_0^{\tau_k} ds \, \phi_s(\lambda_s)) \right) \tag{3.13}$$

$$= \frac{1}{v(x)} E_x^{2(u+g)} \left( v(\lambda_{\tau_k}) \exp(- \int_0^{\tau_k} ds \, 2v(\lambda_s)) \exp(- \int_0^{\tau_k} ds \, \phi_s(\lambda_s)) \right) \tag{3.14}$$

$$= \frac{1}{v(x)} E_x \left( v(B_{\tau_k}) \exp(- \int_0^{\tau_k} ds \, 2(u+g)(B_s)) \exp(- \int_0^{\tau_k} ds \, 2v(B_s)) \right.$$
$$\left. \times \exp(- \int_0^{\tau_k} ds \, \phi_s(\lambda_s)) \right) \tag{3.15}$$

$$= \frac{1}{v(x)} E_x \left( v(B_{\tau_k}) \exp(- \int_0^{\tau_k} ds \, 4g(B_s)) \exp(- \int_0^{\tau_k} ds \, \phi_s(\lambda_s)) \right)$$

$$= \mathbb{Q}_x(\exp - \int_0^\infty dt \, \langle \Upsilon_t^k, \phi_t \rangle; \Upsilon^k \sim 1) \tag{3.16}$$

Line (3.13) comes from the definition of branching at rate $2v$, line (3.14) is from the definition of $\lambda$ as a $v$-process. Killing at rate $2(u+g)$ yields (3.15), and (3.16) follows as in the proof of Theorem 3.2.

By considering the time $\zeta$ of the first branch, we have in a similar manner

$$A_x(n) = \hat{Q}_x \left( \exp - \int_0^\infty dt \, \langle \Upsilon_t^k, \phi_t \rangle; \Upsilon^k \sim n \right)$$

$$= \sum_{j=1}^{n-1} E_x^{2(u+g),v} \left( \exp \left( - \int_0^\zeta ds \, \phi_s(\lambda_s) \right) A_{\lambda(\zeta)}(j) A_{\lambda(\zeta)}(n-j); \zeta < \tau_k \right)$$

$$
\begin{aligned}
&= \sum_{j=1}^{n-1} E_x^{2(u+g),v} \Big( \int_0^{\tau_k} dt\, 2v(\lambda_t) \exp\Big(-\int_0^t ds\, 2v(\lambda_s)\Big) \exp\Big(-\int_0^t ds\, \phi_s(\lambda_s)\Big) \\
&\qquad \times A_{\lambda(t)}(j) A_{\lambda(t)}(n-j) \Big) \\
&= \sum_{j=1}^{n-1} \frac{1}{v(x)} E_x^{2(u+g)} \Big( \int_0^{\tau_k} dt\, 2v^2(\lambda_t) \exp\Big(-\int_0^t ds\, 2v(\lambda_s)\Big) \exp\Big(-\int_0^t ds\, \phi_s(\lambda_s)\Big) \\
&\qquad \times A_{\lambda(t)}(j) A_{\lambda(t)}(n-j) \Big) \\
&= \sum_{j=1}^{n-1} \frac{1}{v(x)} E_x \Big( \int_0^{\tau_k} dt\, 2v^2(\lambda_t) \exp\Big(-\int_0^t ds\, 2v(B_s)\Big) \exp\Big(-\int_0^t ds\, 2(u+g)(B_s)\Big) \\
&\qquad \times \exp\Big(-\int_0^t ds\, \phi_s(\lambda_s)\Big) A_{B(t)}(j) A_{B(t)}(n-j) \Big) \\
&= \sum_{j=1}^{n-1} \frac{1}{v(x)} E_x \Big( \int_0^{\tau_k} dt\, 2v^2(\lambda_t) \exp\Big(-\int_0^t ds\, 4g(B_s)\Big) \exp\Big(-\int_0^t ds\, \phi_s(\lambda_s)\Big) \\
&\qquad \times A_{B(t)}(j) A_{B(t)}(n-j) \Big) \\
&= \mathbb{Q}_x(\exp -\int_0^\infty dt\, \langle \Upsilon_t^k, \phi_t \rangle;\, \Upsilon^k \sim n)
\end{aligned}
\tag{3.18}
$$

Line (3.17) uses the density for the first branch time. Line (3.18) uses induction and the calculation from the proof of Theorem 3.2. $\square$

We now consider a number of examples of such conditioning.

*Example* 3.7. Let $D \subset \mathbb{R}^2$ be a bounded $C^2$-domain. It is shown in [8] that points are hit with positive probability. Thus, when $n = 1$ and $d = 2$, the analogue of the transforms of SV1 would be a conditioning on the event that $\mathcal{R}^D \cap \{z\} \neq \emptyset$. Arguing as in Theorem 5.6 of SV1, we have that

$$
\mathbb{N}_x(e_\phi^k \mid \mathcal{R}^D \cap \{z\} \neq \emptyset) = \frac{1}{g_z(x)} \mathbb{N}_x(e_\phi^k(1 - \exp -\langle X^k, g_z \rangle)),
$$

where $g_z(x) = \mathbb{N}_x(\mathcal{R}^D \cap \{z\} \neq \emptyset)$ satisfies $\Delta u = 4u^2$ in $D$ and has boundary value 0 away from $z$. In fact, $g_z$ is the maximal such solution. Applying Theorem 3.2 or Theorem 3.5, with $u = 0$ and $g = v = g_z$, gives a representation of the conditioned process in terms of a branching tree, throwing off mass which is killed off at rate $g_z$. All branches of the tree converge to $z$.

In view of our earlier discussion, concerning whether $g$ is an $L_{4g}$-potential or has a non-zero $L_{4g}$-harmonic component, it is worth noting the following

**Proposition 3.8.** *For $D \subset \mathbb{R}^2$ a $C^2$ domain, and $z \in \partial D$, the associated tree $\Upsilon$ has infinitely many branches almost surely. In consequence, $g_z$ is an $L_{4g_z}$-potential.*

*Proof.* Remark 1.2 of [11] gives the estimate

$$
\frac{g_z(x)}{d(x,z)^{-2}} \leq M
$$

for every $x \in D$, with a reverse inequality also holding provided $x \in D \cap C$. Here $C \subset D$ is any non-tangential cone with vertex at $z$.

Let $D_k = \{x \mid d(x, \partial D) > 2^{-k}\}$, and $\tau_k = \tau_{D_k}$. It follows that
$$\frac{g_z(x)}{g_z(y)} \leq K, \tag{3.19}$$
as long as $x, y \in C \cap (D_k \setminus D_{k-1})$. A similar statement is true for $K(\cdot, z)$.

Let $\xi_t$ be a $g_z$-transform of the process with generator $L_{2g_z}$, which creates branches at rate $2g_z$. Let $A_k$ be the event that $\xi$ branches between $\tau_{k-1}$ and $\tau_k$. Let $x \in C$ satisfy $d(x, \partial D) = 2^{-k+1}$. Then by (3.19),

$$\begin{aligned}
P_x^{2g_z, g_z}(A_k) &= \frac{1}{g_z(x)} P_x^{2g_z}(g_z(\xi_{\tau_k}), A_k \cap \{\zeta > \tau_k\}) \\
&= \frac{1}{g_z(x)} P_x^{2g_z}\left(g_z(\xi_{\tau_k})\left(1 - \exp -\int_0^{\tau_k} 2g_z(\xi_t)\,dt\right), \zeta > \tau_k\right) \\
&= \frac{1}{g_z(x)} P_x\left(g_z(\xi_{\tau_k})\left(1 - \exp -\int_0^{\tau_k} 2g_z(\xi_t)\,dt\right) \exp -\int_0^{\tau_k} 2g_z(\xi_t)\,dt\right) \\
&= \frac{K(x,z)}{g_z(x)} P_x^{K(\cdot,z)}\left(\frac{g_z(\xi_{\tau_k})}{K(\xi_{\tau_k}, z)}\left(1 - \exp -\int_0^{\tau_k} 2g_z(\xi_t)\,dt\right) \exp -\int_0^{\tau_k} 2g_z(\xi_t)\,dt\right) \\
&\geq P_x^{K(\cdot,z)}\left(\frac{g_z(\xi_{\tau_k})}{g_z(x)} \cdot \frac{K(x,z)}{K(\xi_{\tau_k}, z)}\left(1 - \exp -\int_0^{\tau_k} 2g_z(\xi_t)\,dt\right) \right. \\
&\qquad \left. \times \exp(-\int_0^{\tau_k} 2g_z(\xi_t)\,dt), \xi_{\tau_k} \in C\right) \\
&\geq K P_x^{K(\cdot,z)}\left(\left(1 - \exp -\int_0^{\tau_k} 2g_z(\xi_t)\,dt\right) \exp(-\int_0^{\tau_k} 2g_z(\xi_t)\,dt), \xi_{\tau_k} \in C\right) \\
&\geq K > 0.
\end{aligned}$$

To see the last line, observe that by the Brownian scaling, there is probability $\geq K > 0$ that $\tau_k \in [2^{-2k}, 2^{-2k+1}]$. In combination with the estimates on $g_z$, this shows that there is probability $\geq K > 0$ that $\int_0^{\tau_k} 2g_z(\xi_t)\,dt \in [M^{-1}, M]$.

The statement of the Proposition now follows immediately. □

*Example* 3.9. Consider, more generally, a domain $D \subset \mathbb{R}^d$, and a closed non-polar subset $\Gamma$ of $\partial D$. Let $u = 0$ and take $v(x) = g(x) = \mathbb{N}_x(\mathcal{R}^D \cap \Gamma \neq \emptyset)$. As in the previous example, conditioning the range of super-Brownian motion to hit $\Gamma$ yields the measure $\hat{\mathbb{M}}_x$. Thus Theorem 3.2 or Theorem 3.5 represents this conditioned process in terms of a branching tree $\Upsilon$. Since the mass thrown off by $\Upsilon$ will die before reaching $\Gamma$, we are entitled to interpret $\Upsilon$ as the historical tree of all those "particles" that survive to hit $\Gamma$. We conjecture that, as in Example 3.7, this tree will always have infinitely many branches.

*Example* 3.10. Let $D \subset \mathbb{R}^d$ be a bounded domain, and let $f \geq 0$ be a continuous function on $\partial D$. Consider the solution to $\Delta g = 4g^2$ on $D$, given by
$$g(x) = \mathbb{N}_x(1 - \exp -\langle X^D, f\rangle)$$
(if $\partial D$ is regular, then $g$ has boundary value $f$; see Dynkin [5] or Le Gall [9]). Let $u = 0$ and $v = g$, and consider the associated transform $\hat{\mathbb{M}}_x$. Theorem 3.2 or Theorem 3.5 gives a representation of the solution in terms of a branching backbone $\Upsilon$. Because $g$ is bounded, the rate at which $\Upsilon$ branches is also bounded, and so it will have finitely many branches a.s. Thus $g$, considered as an $L_{4g}$-superharmonic function, will have a non-zero $L_{4g}$-harmonic component.

More generally, we conjecture that $\Upsilon$ will have finitely many branches whenever $g \geq 0$ solves $\Delta g = 4g^2$ and satisfies the condition of Dynkin, that it be dominated by a harmonic (that is, $L_0$-harmonic) function (see [10]).

*Example* 3.11. Let Let $D \subset \mathbb{R}^2$ be a bounded $C^2$-domain. In [11], Le Gall classifies all solutions to $\Delta g = 4g^2$. They are in one-to-one correspondence with pairs $(\Gamma, \nu)$, where $\Gamma$ is a closed subset of $\partial D$, and $\nu$ is a Radon measure on $\partial D \setminus \Gamma$. In the case that $\nu(dz) = f(z)\sigma(dz)$, where $\sigma$ is surface area, the representation of solutions takes the form

$$g(x) = \mathbb{N}_x(\mathcal{R}^D \cap \Gamma \neq \emptyset) + \mathbb{N}_x(1 - \exp -\langle X^D, f \rangle, \mathcal{R}^D \cap \Gamma = \emptyset).$$

Let $u = \mathbb{N}_x(\mathcal{R}^D \cap \Gamma \neq \emptyset)$, so $\frac{1}{2}\Delta u = 2u^2$, and set $v = g - u$. For $\Phi \in \mathcal{F}_k$, the $g$-transformed measure

$$\check{\mathbb{N}}_x^g(\Phi) = \frac{1}{g(x)} \mathbb{N}_x(\Phi(1 - e_g^k))$$

becomes a superposition

$$\check{\mathbb{N}}_x^g = \frac{u(x)}{g(x)} \check{\mathbb{N}}_x^u + \frac{v(x)}{g(x)} \check{\mathbb{N}}_x^v,$$

where

$$\check{\mathbb{N}}_x^u(\Phi) = \frac{1}{u(x)} \mathbb{N}_x(\Phi(1 - e_u^k))$$

$$\check{\mathbb{N}}_x^v(\Phi) = \frac{1}{v(x)} \mathbb{N}_x(\Phi e_u^k (1 - e_v^k)).$$

The measure $\check{\mathbb{N}}_x^u$ is of the type considered in Example 3.9, and is represented in terms of a tree whose branches (conjecturally infinitely many) terminate in $\Gamma$, and throw off mass which is killed at rate $u$.

On the other hand, $\check{\mathbb{N}}_x^v$ is also of the form (3.1) (with $g$, $u$, $v$ all as described above). Thus Theorem 3.2 or Theorem 3.5 give a representation in terms of a tree throwing off mass which gets killed at rate $g$. Each branch of the tree follows a $v$-transform of the process with generator $L_{2(u+g)}$, with branching at rate $2v$. If $v$ is bounded (for example, if $f$ is bounded and vanishes on a neighborhood of $\Gamma$), it then follows immediately that the tree has only finitely many branches, none of which terminate in $\Gamma$. As in the preceding example, we conjecture that the latter property is generic.

*Remark* 3.12. It is natural to ask for relationships between the transforms $\hat{\mathbb{M}}_x$ of (3.1), and $\mathbb{M}_x$ of SV1.

The martingale $\mathbb{M}_x$ from SV1 is defined by the following. Suppose we have $n$ positive solutions in $D$ to the linear equation $L_{4g}v = 0$, labeled $v^1, \ldots, v^n$. Recall that $U^{4g}$ is the potential operator for the generator $L_{4g}$ in $D$, and recursively define a family of functions $v^A$, for $\emptyset \neq A \subseteq N = \{1, \ldots n\}$, as follows:

$$v^A = \begin{cases} v^i & A = \{i\}, \\ 2\sum_{\substack{B \subseteq A \\ \emptyset, A \neq B}} U^{4g}(v^B v^{A \setminus B}) & |A| \geq 2. \end{cases} \quad (3.20)$$

Note that $v^A$ is either finite everywhere on $D$, or $v^A \equiv \infty$.

Let $D_k \Uparrow D$ be an increasing sequence of bounded, smooth subdomains and for each $k$ let $X^k$ be the exit measure from $D_k$, and $e_\phi^k = \exp-\langle X^k, \phi \rangle$. Let $\tau_k$ denote the exit time of a path from $D_k$. For $\emptyset \neq A \subseteq N$ we define

$$M_k^A = \sum_{\sigma \in \mathcal{P}(A)} \exp(-\langle X^k, g \rangle) \prod_{C \in \sigma} \langle X^k, v^C \rangle. \tag{3.21}$$

Going in one direction, we can set $f_\epsilon = \epsilon f$, and let $\epsilon \to 0$. Not only should we recover the measure $\mathbb{M}_x$ (with $n = 1$ and with $v^N$ the harmonic function with boundary value $f$) as a limit of the resulting $\hat{\mathbb{M}}_x$'s, but the probabilistic representation carries over as well, as the branching rate goes to zero leaving a single particle in the limit.

Alternatively, in view of the discussion of Remark 3.4, it should be possible to recover the above martingale $\hat{M}^k$ from martingales of the type $M^k$. Let $h$ be the $L_{4g}$-harmonic function with boundary value $f$. Set $v^i = a_i h$, $i = 1, \ldots, n$, and then define $M^k$ as in (3.21). We conjecture that a suitable choice of constants $a_n$ will ensure that $M^k \to \hat{M}^k$ as $n \to \infty$.

## 4. The martingale $\check{M}_k$.

Finally, we will define another transform, combining features of both the transform $\mathbb{M}_x$ of SV1, and the $\hat{\mathbb{M}}_x$ of (3.1).

Let $D \subset \mathbb{R}^d$ be a domain. Let $n \geq 1$, and suppose that for every nonempty $A \subset N = \{1, \ldots, n\}$, we are given a solution $u^A \geq 0$ to the equation $\frac{1}{2}\Delta u = 2u^2$. Define

$$v_A = \sum_{\substack{N \setminus A \subseteq B \subseteq N \\ B \neq \emptyset}} (-1)^{|A|+|B|+n+1} u^B.$$

Suppose also that the relations

$$v_A \geq 0 \tag{4.1}$$

hold for every $\emptyset \neq A \subseteq N$. Then

**Lemma 4.1.**

$$(a) \quad u^A = \sum_{\substack{B \subseteq N \\ A \cap \overline{B} \neq \emptyset}} v_B,$$

$$(b) \quad \frac{1}{2}\Delta v_A = 4 u^N v_A - 2 \sum_{\substack{B \cup C = A \\ B, C \neq \emptyset}} v_B v_C.$$

*Proof.* To show part (a), observe that

$$\sum_{\substack{B \subseteq N \\ A \cap \overline{B} \neq \emptyset}} v_B = \sum_{\substack{B \subseteq N \\ A \cap \overline{B} \neq \emptyset}} \sum_{\substack{N \setminus B \subseteq C \subseteq N \\ C \neq \emptyset}} (-1)^{|B|+|C|+n+1} u^C$$

$$= \sum_{\emptyset \neq C \subseteq N} (-1)^{|C|+n+1} u^C \sum_{\substack{N \setminus C \subseteq B \subseteq N \\ B \cap A \neq \emptyset}} (-1)^{|B|}$$

$$= \sum_{\emptyset \neq C \subseteq N} (-1)^{|C|+n+1} u^C \left( \sum_{N \setminus C \subseteq B \subseteq N} (-1)^{|B|} - \sum_{N \setminus C \subseteq B \subseteq N \setminus A} (-1)^{|B|} \right) \tag{4.2}$$

$$= \sum_{\emptyset \neq C \subseteq N} (-1)^{|C|+n+1} u^C (-1)^{n-|A|+1} \mathbf{1}_{A=C} \tag{4.3}$$

$$= u^A,$$

where (4.2) follows by rewriting the previous summation and (4.3) holds by virtue of Lemma 2.1.

Thus,

$$\frac{1}{2}\Delta v_A = \sum_{\substack{N\setminus A \subseteq B \subseteq N \\ B \neq \emptyset}} (-1)^{|A|+|B|+n+1} 2(u^B)^2$$

$$= 2 \sum_{\substack{N\setminus A \subseteq B \subseteq N \\ B \neq \emptyset}} (-1)^{|A|+|B|+n+1} \sum_{\substack{C,C' \subseteq N \\ B \cap C, B \cap C' \neq \emptyset}} v_C v_{C'}$$

$$= 2 \sum_{\substack{C,C' \subseteq N \\ C,C' \neq \emptyset}} v_C v_{C'} (-1)^{|A|+n+1} \sum_{\substack{N\setminus A \subseteq B \subseteq N \\ B \cap C, B \cap C' \neq \emptyset}} (-1)^{|B|}$$

$$= 2 \sum_{\substack{C,C' \subseteq N \\ C,C' \neq \emptyset}} v_C v_{C'} (-1)^{|A|+n+1} \Bigg( \sum_{N\setminus A \subseteq B \subseteq N} (-1)^{|B|} - \sum_{N\setminus A \subseteq B \subseteq N\setminus C} (-1)^{|B|}$$

$$- \sum_{N\setminus A \subseteq B \subseteq N\setminus C'} (-1)^{|B|} + \sum_{N\setminus A \subseteq B \subseteq N\setminus (C \cup C')} (-1)^{|B|} \Bigg)$$

$$= 2 \sum_{\substack{C,C' \subseteq N \\ C,C' \neq \emptyset}} v_C v_{C'} (-1)^{|A|+n+1} \Big( -\mathbf{1}_{C=A} - \mathbf{1}_{C'=A} + \mathbf{1}_{C \cup C' = A} \Big) (-1)^{n-|A|}$$

$$= 4 v_A \Bigg( \sum_{\emptyset \neq C \subseteq N} v_C \Bigg) - 2 \sum_{\substack{C \cup C' = A \\ C,C' \neq \emptyset}} v_C v_{C'}$$

$$= 4 u^N v_A - 2 \sum_{\substack{C \cup C' = A \\ C,C' \neq \emptyset}} v_C v_{C'}$$

□

*Remark* 4.2. Though we will not need it, the analogue of the $v^A$ of (3.20) are really

$$v^A = \sum_{\emptyset \neq B \subseteq A} (-1)^{|B|+1} u^B.$$

The following relations could be proved just as above:

$$v^A = \sum_{A \subseteq B \subseteq N} v_B$$

$$v_A = \sum_{A \subseteq B \subseteq N} (-1)^{|A|+|B|} v^B$$

$$u^A = \sum_{\emptyset \neq B \subseteq A} (-1)^{|B|+1} v^B.$$

Now set
$$\check{M}_k = 1 + \sum_{\emptyset \neq A \subseteq N} (-1)^{|A|} \exp -\langle X^k, u^B \rangle.$$

It follows immediately that $\check{M}_k$ is a $\mathbb{N}_x$-martingale, and so for $\Phi \in \mathcal{F}_k$ we can define
$$\check{\mathbb{M}}_x(\Phi) = \frac{1}{v_N(x)} \mathbb{N}_x(\Phi \check{M}_k).$$

**Lemma 4.3.**
$$\check{M}_k = \exp(-\langle X^k, u^N \rangle) \sum_{m=1}^{\infty} \frac{1}{m!} \sum_{\substack{C_1 \cup \cdots \cup C_m = N \\ C_i \neq \emptyset \, \forall i}} \prod_{i=1}^{m} \langle X^k, v_{C_i} \rangle.$$

*Proof.* Write $u^\emptyset = 0$. Then by Lemma 2.2
$$\check{M}_k = \sum_{B \subseteq N} (-1)^{|B|} \exp -\langle X^k, u^B \rangle$$
$$= e_{u^N}^k \sum_{B \subseteq N} (-1)^{|B|} \exp\langle X^k, u^N - u^B \rangle$$
$$= e_{u^N}^k \sum_{B \subseteq N} \sum_{m=0}^{\infty} \frac{(-1)^{|B|}}{m!} \langle X^k, u^N - u^B \rangle^m$$
$$= e_{u^N}^k \left( \left( \sum_{B \subseteq N} (-1)^{|B|} \right) + \sum_{m=1}^{\infty} \sum_{B \subseteq N} \frac{(-1)^{|B|}}{m!} \left( \sum_{\emptyset \neq C \subseteq N \setminus B} \langle X^k, v_C \rangle \right)^m \right)$$
$$= e_{u^N}^k \sum_{m=1}^{\infty} \sum_{B \subseteq N} \frac{(-1)^{|B|}}{m!} \sum_{\substack{C_1, \ldots, C_m \subseteq N \setminus B \\ \emptyset \neq C_i \, \forall i}} \prod_{i=1}^{m} \langle X^k, v_{C_i} \rangle$$
$$= e_{u^N}^k \sum_{m=1}^{\infty} \sum_{\substack{C_1, \ldots, C_m \subseteq N \\ \emptyset \neq C_i \, \forall i}} \frac{1}{m!} \prod_{i=1}^{m} \langle X^k, v_{C_i} \rangle \sum_{B \subseteq N \setminus \cup C_i} (-1)^{|B|}$$
$$= e_{u^N}^k \sum_{m=1}^{\infty} \sum_{\substack{C_1 \cup \cdots \cup C_m = N \\ \emptyset \neq C_i \, \forall i}} \frac{1}{m!} \prod_{i=1}^{m} \langle X^k, v_{C_i} \rangle.$$

$\square$

To describe the probabilistic representation of $\check{\mathbb{M}}_x$, we construct a measure $\tilde{\mathbb{N}}_x$, as before. It has a tree backbone $\Upsilon$, and throws off mass which gets killed at rate $u^N$. In other words, we use $\tilde{\mathbb{N}}_x(\Phi_k) = \mathbb{N}_x(\Phi_k e_{u^N}^k)$, for $\Phi_k \in \mathcal{F}_k$. To construct the backbone, we start a single particle off at $x$, following a $v_N$-transform of the process with generator $L_{4u^N}$. When it dies, say at a point $y$, we choose a pair $(A, A')$ such that $A \cup A' = N$ and $A, A' \neq \emptyset$, according to the law
$$p(A, A'; N)(y) = \frac{v_A(y) v_{A'}(y)}{\sum_{\substack{B \cup B' = N \\ B, B' \neq \emptyset}} v_B(y) v_{B'}(y)}.$$

At its death, the $v^N$-particle splits into a $v^A$-particle and a $v^{A'}$ particle. The $v^A$-particle follows a $v_A$-transform of $L_{4u^N}$ and when it dies, it splits into a $v^B$-particle and a $v^{B'}$-particle, where $(B, B')$ is chosen according to law $p(B, B'; A)$, and so on.

This gives us a tree $\Upsilon$ of branching particles, each tagged with a set $A$. We may form $\Upsilon^k$ as before, by pruning off all particles (together with their descendants), once they leave $D_k$. We write $\check{Q}_x$ for the law of $\Upsilon$, and set

$$\check{\mathbb{N}}(\exp-\langle Y^k, \phi \rangle) = \check{Q}_x(\exp - \int_0^\infty dt\, 4\langle \Upsilon_t^k, \check{\mathbb{N}}.(1 - e_\phi^k)\rangle).$$

**Theorem 4.4.** *Assume condition* (4.1). *Then*

$$\check{\mathbb{M}}_x\left(\exp-\langle X^k, \phi\rangle\right) = \check{\mathbb{N}}_x\left(\exp-\langle Y^k, \phi\rangle\right)$$

*Remark* 4.5. Using historical processes, as in the last section of SV1 one can show that $\check{\mathbb{M}}_x = \check{\mathbb{N}}_x$ on $\mathcal{F}_k$.

*Proof.* In the present context, it is useful to label all the particles of $\Upsilon^k$ that exit $D_k$, by placing an order on them. So let $F_k$ be the set of such particles, and set $\gamma_k = |F_k|$. For $A \subseteq N$, let

$$\mathcal{S}_m(A) = \{(C_1, \ldots, C_m) : C_1 \cup \cdots \cup C_m = A, \emptyset \neq C_i\ \forall i\}.$$

If $\gamma_k = m$, choose at random an ordering of $F_k$, and for $\Lambda = (C_1, \ldots, C_m) \in \mathcal{S}_m(N)$, write $\Upsilon^k \approx \Lambda$ for the event that the $i$th particle is tagged with the set $C_i$, $i = 1, \ldots, m$. Thus for example,

$$\check{Q}_x(\gamma_k = m) = \sum_{\Lambda \in \mathcal{S}_m(N)} \check{Q}_x(\Upsilon^k \approx \Lambda). \tag{4.4}$$

Note that if $M \subseteq S$ are sets with $|S| = m$ and $|M| = j$, then there are $\binom{m}{j}$ orderings of $S$ compatible with any given orders on $M$ and on $S \setminus M$. In other words, if $\sigma$ is any order on $S$, and if $\Sigma$ is an order on $S$ picked at random, then the conditional probability

$$P(\Sigma = \sigma \mid \Sigma_M = \sigma_M, \Sigma_{S \setminus M} = \sigma_{S \setminus M}) = 1/\binom{m}{j} \tag{4.5}$$

(writing $\sigma_M$ etc ... for the restriction of $\sigma$ to $M$).

As described initially, the root particle of the tree is always a $v^N$-particle. It is convenient, for purposes of induction, to allow the same notation to cover the situation that we start with our root being a $v^A$-particle for some $A \subseteq N$. In this case, (4.4) still holds, but with $\Lambda \in \mathcal{S}_m(N)$ replaced by $\Lambda \in \mathcal{S}_m(A)$. With this in mind, we may define another restriction operation as follows. For $1 \leq i_1 < \cdots < i_k \leq m$, set

$$(C_1, \ldots, C_m)|_{\{i_1, \ldots, i_k\}} = (C_{i_1}, \ldots, C_{i_k}).$$

Thus, if $\Lambda = (C_1, \ldots, C_m) \in \mathcal{S}_m(A)$ and $M \subseteq \{1, \ldots, m\}$, we will have that $\Lambda|_M \in \mathcal{S}_m(B)$, for $B = \cup_{i \in M} C_i$. As a shorthand for the latter, we write $\Lambda(M) = \cup_{i \in M} C_i$.

We will show, by induction on $m \geq 1$, that for $\emptyset \neq A \subseteq N$, and $(C_1, \ldots, C_m) \in \mathcal{S}_m(A)$,

$$\check{Q}_x(\exp - \int_0^\infty dt\, 4\langle \Upsilon_t^k, \check{\mathbb{N}}.(1 - e_\phi^k)\rangle; \Upsilon^k \approx (C_1, \ldots, C_m))$$

$$= \frac{1}{m!\, v_A(x)} \mathbb{N}_x(e^k_{\phi + u^N} \prod_{i=1}^m \langle X^k, v_{C_i}\rangle). \tag{4.6}$$

Taking $A = N$ and summing over $\mathcal{S}_m(N)$ will then establish the theorem.

The initial stage of the induction, with $m = 1$ follows exactly as in the proof of Theorem 3.2. So let $m > 1$ and assume the inductive hypothesis for all $A \subseteq N$, and for all values smaller than $m$. For simplicity, we will verify (4.6) in the case $A = N$. For $\zeta$ the lifetime of the initial particle, and $\Lambda = (C_1, \ldots, C_m)$, we have that

$$\check{\mathbb{Q}}_x(\exp - \int_0^\infty dt\, 4\langle \Upsilon_t^k, \tilde{\mathbb{N}}.(1 - e_\phi^k)\rangle, \Upsilon^k \approx \Lambda)$$

$$= \sum_{j=1}^{m-1} \sum_{\substack{M \subset \{1,\ldots,m\} \\ |M|=j}} E_x^{4u^N, v_N}\left(1_{\zeta < \tau_k} \exp\left(-\int_0^\zeta dt\, 4\tilde{\mathbb{N}}_{\xi_\zeta}(1 - e_\phi^k)\right)\right.$$

$$\times p(\Lambda(M), \Lambda(M^c); N)(\xi_\zeta) \binom{m}{j}^{-1} \tag{4.7}$$

$$\times \check{Q}_{\xi_\zeta}\left(\exp - \int_0^\infty dt\, 4\langle \Upsilon_t^k, \tilde{\mathbb{N}}.(1 - e_\phi^k)\rangle, \Upsilon^k \approx \Lambda|_M\right)$$

$$\left.\times \check{Q}_{\xi_\zeta}\left(\exp - \int_0^\infty dt\, 4\langle \Upsilon_t^k, \tilde{\mathbb{N}}.(1 - e_\phi^k)\rangle, \Upsilon^k \approx \Lambda|_{\{1,\ldots,m\}\setminus M}\right)\right)$$

$$= \sum_{j=1}^{m-1} \sum_{\substack{M \subset \{1,\ldots,m\} \\ |M|=j}} \frac{1}{v_N(x)} E_x^{4u^N}\left(\int_0^{\tau_k} ds\left(2 \sum_{(A,A') \in \mathcal{S}_2(N)} v_A(\xi_s) v_{A'}(\xi_s)\right) 1_{\zeta > s}\right.$$

$$\times \frac{j!(m-j)!}{m!} p(\Lambda(M), \Lambda(M^c); N)(\xi_s)$$

$$\times \exp\left(-\int_0^s dt\, 4(\tilde{\mathbb{N}}_{\xi_\zeta}(1 - e_{\phi+u^N}^k) - \tilde{\mathbb{N}}_{\xi_\zeta}(1 - e_{u^N}^k))\right) \tag{4.8}$$

$$\times \check{Q}_{\xi_s}\left(\exp - \int_0^\infty dt\, 4\langle \Upsilon_t^k, \tilde{\mathbb{N}}.(1 - e_\phi^k)\rangle, \Upsilon^k \approx \Lambda|_M\right)$$

$$\left.\times \check{Q}_{\xi_s}\left(\exp - \int_0^\infty dt\, 4\langle \Upsilon_t^k, \tilde{\mathbb{N}}.(1 - e_\phi^k)\rangle, \Upsilon^k \approx \Lambda|_{\{1,\ldots,m\}\setminus M}\right)\right)$$

$$= \sum_{j=1}^{m-1} \sum_{\substack{M \subset \{1,\ldots,m\} \\ |M|=j}} \frac{2}{m! v_N(x)} E_x\left(\int_0^{\tau_k} ds\, \exp\left(-\int_0^s dt\, 4u^N(\xi_t)\right) \exp\left(\int_0^s dt\, 4u^N(\xi_t)\right)\right.$$

$$\times \mathcal{N}_s(e_{\phi+u^N}^k) j! v_{\Lambda(M)}(\xi_s) \check{Q}_{\xi_s}\left(\exp - \int_0^\infty dt\, 4\langle \Upsilon_t^k, \tilde{\mathbb{N}}.(1 - e_\phi^k)\rangle, \Upsilon^k \approx \Gamma|_M\right)$$

$$\left.\times (m-j)! v_{\Lambda(M^c)}(\xi_s) \check{Q}_{\xi_s}\left(\exp - \int_0^\infty dt\, 4\langle \Upsilon_t^k, \tilde{\mathbb{N}}.(1 - e_\phi^k)\rangle, \Upsilon^k \approx \Gamma|_{M^c}\right)\right) \tag{4.9}$$

$$= \sum_{\substack{M \subset \{1,\ldots,m\} \\ 1 < |M| < m}} \frac{2}{m! v_N(x)} E_x\left(\int_0^{\tau_k} ds\, \mathcal{N}_s(e_{\phi+u^N}^k) \mathbb{N}_{\xi_s}\left(e_{\phi+u^N}^k \prod_{i \in M} \langle X^k, v_{C_i}\rangle\right)\right.$$

$$\left.\times \mathbb{N}_{\xi_s}\left(e_{\phi+u^N}^k \prod_{i \in \{1,\ldots,m\}\setminus M} \langle X^k, v_{C_i}\rangle\right)\right) \tag{4.10}$$

$$= \frac{1}{m! v_N(x)} \mathbb{N}_x (e^k_{\phi+u^N} \prod_{i=1}^{m} \langle X^k, v_{C_i} \rangle). \tag{4.11}$$

Here, line (4.7) follows from (4.5) and the definition of $\check{Q}$, (4.8) follows from (2.4), (4.9) from (2.1) and Lemma 2.3, (4.10) by the inductive hypothesis, and (4.11) from (3.3). □

*Example* 4.6. Let $\Gamma_1, \ldots, \Gamma_n$ be disjoint closed non-polar subsets of $\partial D$, where $D \subset \mathbb{R}^d$ is smoothly bounded. Set

$$u^A(x) = \mathbb{N}_x(\mathcal{R}^D \cap \bigcup_{i \in A} \Gamma_i \neq \emptyset).$$

Then, by inclusion-exclusion,

$$v_A(x) = \mathbb{N}_x(\mathcal{R}^D \cap \Gamma_i \neq \emptyset \ \forall i \in A, \mathcal{R}^D \cap \Gamma_i = \emptyset \ \forall i \in N \setminus A).$$

Thus (4.1) holds (and, referring to Remark 4.2, $v^A(x) = \mathbb{N}_x(\mathcal{R}^D \cap \Gamma_i \neq \emptyset \ \forall i \in A))$. By (5.20) of SV1 we have that in fact,

$$\check{\mathbb{M}}_x(\Phi_k) = \mathbb{N}_x(\Phi_k \mid \mathcal{R}^D \cap \Gamma_i \neq \emptyset \ \forall i \in N),$$

so that Theorem 4.4 provides a particle representation of the process conditioned to hit each $\Gamma_i$.